\documentclass{amsart}

\usepackage[all]{xy}
\usepackage{amssymb,latexsym}

\newtheorem{theorem}{Theorem}
\newcommand{\bt}{\begin{theorem}}
\newcommand{\et}{\end{theorem}}
\newtheorem{lemma}{Lemma}
\newcommand{\bl}{\begin{lemma}}
\newcommand{\el}{\end{lemma}}
\newtheorem{corollary}{Corollary}
\newcommand{\bc}{\begin{corollary}}
\newcommand{\ec}{\end{corollary}}
\newcommand{\beq}{\begin{equation}}
\newcommand{\eeq}{\end{equation}}
\newcommand{\benum}{\begin{enumerate}}
\newcommand{\eenum}{\end{enumerate}}

\newcommand{\N}{\ensuremath{ \mathbf N }}
\newcommand{\Z}{\ensuremath{\mathbf Z}}

\newcommand{\R}{\ensuremath{\mathbf R}}

\newcommand{\mcb}{\ensuremath{ \mathcal B}}

\newcommand{\bmat}{\left(\begin{matrix}}
\newcommand{\emat}{\end{matrix}\right)}
\DeclareMathOperator{\qand}{\quad\text{and}\quad}
\DeclareMathOperator{\qqand}{\qquad\text{and}\qquad}

\title[Fourth element in the greedy $B_h$-set]{The fourth positive element in the greedy $B_h$-set}
\author{Melvyn B. Nathanson}
\address{Department of Mathematics\\Lehman College (CUNY)\\Bronx, NY 10468} 
\email{melvyn.nathanson@lehman.cuny.edu} 

\author{Kevin O'Bryant}
\address{Department of Mathematics\\College of Staten Island (CUNY)\\ Staten Island, NY 10314} 
\email{obryant@csi.cuny.edu}



\date{\today}

\begin{document}
\maketitle

\begin{abstract}
For $h \geq 1$, a $B_h$-set is a set of integers such that every integer $n$ has at most one representation 
in the form $n = a_{i_1} + \cdots + a_{i_h}$, where $a_{i_r} \in A$ for all $r = 1,\ldots, h$ 
and $a_{i_1} \leq \ldots \leq a_{i_h}$. 
The greedy $B_h$-set is the infinite set of nonnegative integers 
$\{a_0(h), a_1(h), a_2(h), \ldots \}$ 
constructed as follows:  If $a_0(h) = 0$ and $\{a_0(h), a_1(h), a_2(h), \ldots, a_k(h) \}$ 
is a $B_h$-set, then $a_{k+1}(h)$ is the least positive integer such that 
$\{a_0(h), a_1(h), a_2(h), \ldots, a_k(h), a_{k+1}(h) \}$ is a $B_h$-set. 
Then $a_1(h) = 1$, $a_2(h) = h+1$, and $a_3(h) = h^2+h+1$ for all $h$.  
This paper proves that $a_4(h)$, the fourth term of the greedy $B_h$-set is 
$\left( h^3 + 3h^2 + 3h + 1\right) /2$ if $h$ is odd and 
$\left( h^3 + 2h^2 + 3h + 2\right) /2$ if $h$ is even. 
\end{abstract}

\section{The greedy algorithm}
Let $h$ be a positive integer.  
A finite or infinite set $A$  of nonnegative integers is a $B_h$-set if no integer has two different 
representations as sums of $h$ elements of $A$. 
Equivalently, the set $A$ is a $B_h$-set if the equation 
\[
a_{i_1}+  \cdots + a_{i_h} = a_{j_1}+  \cdots + a_{j_h} 
\]
with $a_{i_1},\ldots, a_{i_h}, a_{j_1},\ldots, a_{j_h} \in A$ and 
\[
a_{i_1} \leq \cdots \leq  a_{i_h} 
\]
and 
\[
a_{j_1} \leq \cdots \leq a_{j_h} 
\]
implies $a_{i_r} = a_{j_r}$ for all $r  = 1,\ldots, h$.  
A $B_2$-set is also called a Sidon set. 

For every positive integer $h$, we  use a greedy algorithm to construct a $B_h$-set 
$\{a_k(h): k = 0,1,2,\ldots\}$ as follows:  $a_0(h) = 0$ and, 
if $\{a_0(h), a_1(h), \ldots, a_k(h)\}$ is a $B_h$-set, then $a_{k+1}(h)$ is 
the least positive integer such that  $\{a_0(h), a_1(h), \ldots, a_k(h), a_{k+1}(h)\}$ 
is a $B_h$-set. 
The greedy $B_1$-set is simply the set $\N_0$ of nonnegative integers, that is, $a_k(1) = k$ for all
$k \in \N_0$.  The greedy $B_2$-set is the (shifted) Mian-Chowla sequence 
(Mian-Chowla~\cite{mian-chow44}, Guy~\cite[Section E28]{guy04}).

For all $h \geq 1$ we have 
\[
a_0(h) = 0,\quad a_1(h) = 1, \qand a_2(h) = h+1.  
\]
From computer calculations of initial segments of $B_h$-sets for small $h$, O'Bryant~\cite{obry23a} 
conjectured that, for all positive integers $h$, 
\[
a_3(h) = h^2+h+1. 
\]
This formula for $a_3(h)$ is proved in~\cite{nath2023x}.  
  In this paper we obtain an explicit quasi-polynomial for $a_4(h)$. 
  
  \bt                                                       \label{Bh:theorem} 
The fourth positive element in the greedy $B_h$-set is 
\begin{align*}
a_4(h) 
& = \begin{cases}
\frac{1}{2} \left( h^3 + 3h^2 + 3h + 1\right)  & \text{if $h$ is odd} \\
\frac{1}{2}\left( h^3 + 2h^2 + 3h + 2\right)   & \text{if $h$ is even.}
\end{cases}
\end{align*} 
\et

Using the floor function, we may also write 
\[
a_4(h) = \left[ \frac{h+3}{2} \right] h^2 +  \left[ \frac{3h}{2} \right] +1.
\]
It is an open problem to compute exact values or even asymptotic estimates 
for $a_k(h)$ for $k \geq 5$. 
Indeed, it is not even known if $a_k(h) < a_k(h+1)$ for all integers $h \geq 1$ 
and $k \geq 2$.  We have the upper bound  (in~\cite{nath2023x}) 
\[
a_k(h) \leq \sum_{i=0}^{k-1} h^i < h^{k-1} + 2h^{k-2}
\]
 for all positive integers $h$ and $k$.

\section{Lower bound for $a_4(h)$}
In this section we compute a lower bound for $a_4(h)$.

For $u,v \in \R$,   define the interval of integers $[u,v] = \{ n \in \Z: u \leq n \leq v\}$. 

\bl                                      \label{Bh:lemma:1} 
For all integers $h \geq 2$, 
\[
a_4(h) \geq    \begin{cases}
\frac{1}{2}\left( h^3 + 3h^2 + 3h + 1\right)  & \text{if $h$ is odd}\\
\frac{1}{2}\left( h^3 + 2h^2 + 3h + 2\right)   & \text{if $h$ is even.}
\end{cases}
\]
\el

\begin{proof}

Let  $\mcb$ be the set of positive integers $b$ such that $a_4(h) \neq b$. 
A sufficient condition that  $b \in \mcb$ is the existence 
of nonnegative integers $x_1,x_2,y_1,y_2,y_3$ with  
\[
x_1+x_2  \leq h-1, \qquad y_1+y_2+y_3 \leq h, \qquad x_1y_1 = x_2y_2 = 0
\]
such that 
\[
b + x_1 + x_2(h+1) = y_1 + y_2(h+1) + y_3(h^2+h+1)
\]
or, equivalently, 
\[ 
b =  y_3(h^2+h+1)  + (y_2 - x_2) (h+1) + (y_1-x_1). 
\]
We consider separately the two cases: $x_2 = 0$ and $y_2 = 0$. 

Let $x_2 = 0$ and let $1 \leq y_3 \leq h$ and $0\leq y_2 \leq h-y_3$.    
Because 
\[
0 \leq x_1 \leq h-1 \qqand 0 \leq y_1 \leq h-y_3-y_2
\]
the set \mcb\ contains the interval  
\begin{align}                          
 y_3(h^2+ & h+1) + y_2(h+1) +  \left[ -(h-1), h-y_3-y_2 \right]  \nonumber   \\
& = y_3(h^2+h+1)  +  \left[  y_2(h+1) -h+1, y_2h + h- y_3 \right].        \label{Bh:beqn-y}
\end{align}   
If 
\beq               \label{Bh:beqn-y2y3}
 0 \leq y_2 \leq h-y_3 - 1 
 \eeq
then  \mcb\ also contains the interval 
\begin{align}                          
 y_3(h^2+ & h+1)  +   \left[   (y_2+1) (h+1) -h +1,  (y_2+1) h + h- y_3 \right]    \nonumber   \\
& = y_3(h^2+h+1)  +  \left[  y_2(h+1) + 2, y_2h  + 2h   - y_3 \right].       \label{Bh:beqn-y2}
\end{align} 
Inequality~\eqref{Bh:beqn-y2y3} implies 
\[
y_2(h+1) + 2  \leq y_2h + h- y_3 +1
\]
and so intervals~\eqref{Bh:beqn-y} and~\eqref{Bh:beqn-y2} overlap.  
Therefore, \mcb\ contains the interval 
\begin{align}  
y_3 & (h^2+h+1)  +  \bigcup_{y_2=0}^{h-y_3} \left[  y_2(h+1) -h+1, y_2h + h- y_3 \right]   \nonumber  \\
& = y_3(h^2+h+1) + \left[ -h+1, (h-y_3) h + h- y_3 \right] \nonumber  \\\
& = \left[  y_3(h^2+h+1) -h+1, y_3h^2 + h^2+h \right].      \label{Bh:yint}
\end{align}

Let $y_2 = 0$ and let $1 \leq y_3 \leq h$ and $0 \leq x_2 \leq h-y_3$. 
Note that $h-y_3 \leq h -1$. 
Because 
\[
0 \leq x_1 \leq h-1-x_2 \qqand 0 \leq y_1 \leq h-y_3 
\]
the set \mcb\ contains the interval 
\begin{align}                          
y_3(h^2+ & h+1)  - x_2(h+1) +  \left[ -(h-1-x_2), h -y_3 \right]  \nonumber   \\
& = y_3(h^2+h+1)  +  \left[  -x_2h -h + 1,  - x_2(h+1) + h -  y_3  \right].          \label{Bh:beqn-x}
\end{align} 
If 
\beq               \label{Bh:beqn-x2y3}
0 \leq x_2 \leq h-y_3-1 
 \eeq
then  \mcb\ also contains the interval 
\begin{align}                          
y_3(h^2+ & h+1)  +  \left[  -(x_2+1)h -h + 1,  - (x_2+1)(h+1)+ h - y_3 \right]    \nonumber   \\
& = y_3(h^2+h+1)  +  \left[  -x_2h - 2h + 1,  - x_2(h+1) - y_3 -1 \right].         \label{Bh:beqn-x2}   
\end{align} 
Inequality~\eqref{Bh:beqn-x2y3} implies 
\[
 -x_2h -h  \leq - x_2(h+1)  - y_3 - 1
\]
and so  intervals~\eqref{Bh:beqn-x}  and~\eqref{Bh:beqn-x2}  overlap.  
Therefore, \mcb\ contains the interval 
\begin{align} 
 y_3 (h^2+ & h+1) + \bigcup_{x_2=0}^{h-y_3}  \left[  -x_2h - h + 1,  - x_2(h+1) +  h-y_3 \right]    \nonumber \\
& = y_3(h^2+h+1) +  \left[  -(h-y_3 )h - h + 1,  h-y_3 \right] \nonumber \\
& =  \left[  y_3(h^2+ 2h+1)  -h^2 - h + 1,  y_3(h^2+h) +  h  \right].      \label{Bh:xint}
\end{align}

Intervals~\eqref{Bh:yint} and~\eqref{Bh:xint} overlap because $y_3 \leq h \leq 2h$, 
and so \mcb\ contains the interval 
\beq 
I(y_3) =   \left[  y_3(h^2+ 2h+1) -h^2 - h +1, y_3h^2 + h^2+h \right]
\eeq
for $y_3 \in [1,h]$.  
These intervals move to the right as $y_3$ increases.  
Intervals $(I(y_3)$ and $I(y_3+1)$ overlap if 
\[
 (y_3+1) (h^2+ 2h+1)  -h^2 - h + 1\leq  y_3h^2+ h^2 +  h  + 1
\]
or, equivalently, if  $(2h+1)y_3 \leq h^2-1$ or   
\[
y_3 \leq \left[ \frac{h^2 -1}{2h+1} \right]  = \begin{cases}
\frac{h-1}{2} & \text{if $h$ is odd } \\ 
\frac{h-2}{2} & \text{if $h$ is even.} 
\end{cases}
\]
For odd  $h \geq 3$, the set \mcb\ contains the interval 
\begin{align}                 \label{Bh:odd1}
\bigcup_{y_3=1}^{(h+1)/2} I(y_3) 
& = \bigcup_{y_3=1}^{(h+1)/2}   \left[  y_3(h^2+ 2h+1)  -h^2 - h + 1,  y_3h^2+h^2 +  h  \right] \nonumber  \\ 
& = \left[  h+2,  \left(\frac{h+1}{2} \right) h^2+ h^2 +  h  \right]  \nonumber \\ 
& = \left[  h+2,  \frac{h^3+ 3h^2 + 2h}{2}  \right]. 
\end{align} 
For even $h \geq 2$, the set \mcb\ contains the interval 
\begin{align}                       \label{Bh:even1}
\bigcup_{y_3=1}^{h/2} I(y_3) & =  \left[  y_3(h^2+ 2h+1)  -h^2 - h + 1,  y_3h^2+h^2 +  h  \right]    \nonumber  \\ 
& = \left[  h+2,  \left(\frac{h}{2} \right) h^2+ h^2 +  h  \right]    \nonumber    \\
& =  \left[  h+2,  \frac{h^3+ 2h^2 + 2h}{2}  \right].
\end{align} 

For odd  $h \geq 3$, let 
\[
y_2 = 0, \qquad y_3 = \frac{h+3}{2}, \qquad x_2 = \frac{h+1}{2}. 
\]
The set \mcb\ contains the interval 
\begin{align}              \label{Bh:odd2}
y_3(h^2 & +h+1) - x_2(h+1) + [-(h-1-x_2), h-y_3] \nonumber  \\
& = \left( \frac{h+3}{2}  \right) (h^2+h+1) - \left( \frac{h+1}{2}  \right)(h+1) + \left[ -\frac{h-3}{2}, \frac{h-3}{2} \right]  \nonumber  \\
& =\left[ \frac{h^3 + 3h^2 + h + 5}{2}, \frac{h^3 + 3h^2 + 3h - 1}{2} \right] 
\end{align}
Intervals~\eqref{Bh:odd1} and~\eqref{Bh:odd2}  
overlap and so \mcb\ contains the interval 
\[
\left[  h+2,  \frac{h^3 + 3h^2 + 3h - 1}{2} \right].   
\]
It follows that, for odd  $h \geq 3$, we have the lower bound 
\[
a_4(h) \geq  \frac{h^3 + 3h^2 + 3h - 1}{2} + 1 =  \frac{h^3 + 3h^2 + 3h + 1}{2}.
\]

For even $h \geq 2$, let 
\[
y_2 = 0, \qquad y_3 = \frac{h+2}{2}, \qquad x_2 = \frac{h}{2}. 
\]
The set \mcb\ contains the interval 
\begin{align}                       \label{Bh:even2}
y_3(h^2 & +h+1) - x_2(h+1) + [-(h-1-x_2), h-y_3 ] \nonumber  \\
& = \left( \frac{h+2}{2} \right)  (h^2+h+1) - \left( \frac{h}{2}\right)  (h+1) 
+ \left[ -\frac{h-2}{2}, \frac{h-2}{2} \right]    \nonumber    \\
& =   \left[ \frac{h^3 + 2h^2 + h + 4}{2}, \frac{h^3 + 2h^2 + 3h}{2} \right] 
\end{align}
Intervals~\eqref{Bh:even1} and~\eqref{Bh:even2}  
overlap and  so \mcb\ contains the interval 
\[
\left[  h+2,  \frac{h^3 + 2h^2 + 3h}{2} \right].   
\]
It follows that, for even  $h \geq 2$, we have the lower bound 
\[
a_4(h) \geq  \frac{h^3 + 2h^2 + 3h}{2} + 1 =  \frac{h^3 + 2h^2 + 3h + 2}{2}.
\]
This completes the proof of Lemma~\ref{Bh:lemma:1}. 
\end{proof}

\section{The upper bound for $a_4(h)$}

In this section we compute an upper bound for $a_4(h)$.

\bl                                          \label{Bh:lemma:2} 
For all positive integers $h$, 
\[
a_4(h) \leq    \begin{cases}
\frac{1}{2}\left( h^3 + 3h^2 + 3h + 1\right)  & \text{if $h$ is odd}\\
\frac{1}{2}\left( h^3 + 2h^2 + 3h + 2\right)   & \text{if $h$ is even.}
\end{cases}
\]
\el

\begin{proof}
Let
\[
H =  \begin{cases}
\frac{1}{2}(h^2 + 2h + 1) & \text{if $h$ is odd}\\
\frac{1}{2}(h^2 + h + 2)& \text{if $h$ is even.}
\end{cases}
\]
We must prove that 
\[
a_4(h) \leq (h+1)H. 
\]
This inequality is equivalent to the statement that there do not exist nonnegative integers 
$x_0, x_1, x_2, x_3$ and $y_1, y_2, y_3$ such that 
\begin{align}                     \label{Bh:eqn}
x_0 (h+1)H & + x_1 + x_2(h+1) + x_3(h^2+h+1) \\ 
& = y_1 + y_2(h+1) + y_3(h^2+h+1) \nonumber
\end{align}
with 
\beq                    \label{Bh:eqn-cond1}
x_0 + x_1 + x_2 + x_3 \leq h, \qquad  y_1 + y_2 + y_3 \leq h
\eeq
and 
\beq                    \label{Bh:eqn-cond2}
x_0 \geq 1, \qquad x_1y_1 = x_2y_2 = x_3y_3 = 0.
\eeq

Suppose that  equation~\eqref{Bh:eqn} has a solution satisfying conditions~\eqref{Bh:eqn-cond1} 
and~\eqref{Bh:eqn-cond2}.   Note that 
\[
h^2 +1 < 2H.
\]
If $x_0 \geq 2$, then 
\begin{align*}
(h+1)(h^2+1) & < 2(h+1)H  \leq x_0(h+1)H \\
& \leq  y_1 + y_2(h+1) + y_3(h^2+h+1) \\
& \leq h(h^2+h+1) \\ 
& < (h+1) (h^2 +1)
\end{align*}
which is absurd.  Therefore, 
\[
x_0=1.
\]
If $y_3 = 0$, then 
\[
 \frac{(h+1)(h^2+1)}{2}  < (h+1)H  \leq y_1+y_2(h+1) \leq h(h+1)
\]
and so $(h-1)^2 < 0$, which is absurd.  Therefore, 
\[
y_3 \geq 1 \qqand x_3 = 0.
\]
With $x_0=1$ and $x_3=0$, equation~\eqref{Bh:eqn}  becomes
\begin{align}                 \label{Bh:eqn-2}
x_1+ (H+x_2)(h+1) & = y_1+y_2(h+1)+y_3(h^2+h+1) \\
& =  y_1+ y_3 + (y_2 + y_3h) (h+1). \nonumber
\end{align}
We obtain the congruence 
\[
x_1 \equiv y_1 + y_3 \pmod{h+1}.
\]
If $x_1 = 0$, then $y_1+y_3 \equiv 0 \pmod{h+1}$.  The inequalities 
$y_1 \geq 0$ and $y_3 \geq 1$ imply  
$y_1+y_3 \geq h+1$, which contradicts condition~\eqref{Bh:eqn-cond1}.  
Therefore, 
\[
x_1 \geq 1 \qqand y_1 = 0 
\]
and 
\[
x_1 \equiv y_3 \pmod{h+1}.
\]
Because $1 \leq y_3 \leq h$, if $x_1 \neq y_3$, then $x_1 \geq y_3+h+1 > h$,
which is absurd.  Therefore, 
\[
x_1 = y_3
\]
and equation~\eqref{Bh:eqn-2} becomes, simply, 
\beq                    \label{Bh:eqn-3}
H  + x_2 = y_2 +y_3h  
\eeq 
where $x_2,y_2, y_3$ are nonnegative integers such that 
\[
x_2 + y_3 \leq h-1,\quad y_2+y_3 \leq h, \qqand x_2y_2 = 0.
\]
We consider separately the two cases: $h$ odd and $h$ even. 

Let $h$ be odd.  If $y_3 \leq (h+1)/2$, 
then equation~\eqref{Bh:eqn-3} gives   
\begin{align*}
\frac{h^2+2h+1}{2} &  \leq y_2 + y_3h \leq (h-y_3) + y_3 h \\ 
& = h + y_3(h-1) \leq h + \frac{h^2-1}{2} \\
& = \frac{h^2+2h-1}{2} 
\end{align*}
which is absurd.  It follows that $y_3 > (h+1)2$ and, because $h$ is odd, that 
\[
y_3 \geq \frac{h+3}{2}.
\]
If $y_2 = 0$, then 
 \[
 x_2 = y_3h-H \geq  \frac{h^2+3h}{2}   -  \frac{h^2 + 2h + 1}{2} \geq \frac{h - 1}{2} 
 \]
and so 
\[
h + 1=   \frac{h-1}{2} + \frac{h+3}{2} \leq  x_2 + y_3 \leq h-1
\]
which is absurd.  Therefore, 
\[
y_2 \geq 1 \qqand x_2 = 0. 
\]
Equation~\eqref{Bh:eqn-3} becomes 
\[
H  = y_2 + y_3h
\]
Equivalently,  
\[
h^2 + 2h + 1= 2H =  2y_2 + 2y_3h  \geq 2y_2 + h^2 + 3h  
\]
and so 
\[
1 \geq 2y_2+h 
\]
which is absurd.  
This completes the proof in the odd case.

Let $h$ be even. 
 If $y_3 \leq h/2$, 
then equation~\eqref{Bh:eqn-3} gives   
\begin{align*}
\frac{h^2+h+2}{2} &  \leq y_2 + y_3h \leq (h-y_3) + y_3 h \\ 
& = h + y_3(h-1) \leq h + \frac{h^2-h}{2} \\
& = \frac{h^2+h}{2} 
\end{align*}
which is absurd.  It follows that $y_3 > h/2$ and, because $h$ is even, that 
\[
y_3 \geq \frac{h+2}{2}.
\] 
If $y_2 = 0$, then 
 \[
 x_2 = y_3h- H \geq \frac{h^2+2h}{2} - \frac{h^2 + h + 2}{2} =  \frac{h - 2}{2}
 \]
and so 
\[
h =  \frac{h+2}{2} +   \frac{h-2}{2} \leq  y_3+x_2 \leq h-1
\]
which is absurd.  Therefore, 
\[
y_2 \geq 1 \qqand x_2 = 0.
\]
Equation~\eqref{Bh:eqn-3} becomes 
\[
H  = y_2 + y_3h
\]
Equivalently,  
\[
h^2 + h + 2 = 2H = 2y_2 + 2y_3h \geq 2y_2 + h^2 + 2h   
\]
and so 
\[
2 \geq 2y_2+h 
\]
which is absurd. 
This completes the proof in the even case.  
\end{proof}

Theorem~\ref{Bh:theorem}, the exact formula for $a_4(h)$, 
follows immediately from Lemmas~\ref{Bh:lemma:1} and~\ref{Bh:lemma:2}. 

\def\cprime{$'$} \def\cprime{$'$} \def\cprime{$'$}
\providecommand{\bysame}{\leavevmode\hbox to3em{\hrulefill}\thinspace}
\providecommand{\MR}{\relax\ifhmode\unskip\space\fi MR }
\providecommand{\MRhref}[2]{
  \href{http://www.ams.org/mathscinet-getitem?mr=#1}{#2}
}
\providecommand{\href}[2]{#2}

\end{document}